\documentclass[11pt]{article}
\usepackage[margin=1.15in]{geometry}
\usepackage[T1]{fontenc}
\usepackage{lmodern}
\usepackage{amsmath,amssymb,amsthm,mathtools}
\usepackage{mathrsfs}
\usepackage{enumitem}
\usepackage{microtype}
\IfFileExists{stmaryrd.sty}{\usepackage{stmaryrd}}{%
  \newcommand{\llbracket}{\mathopen{[\mkern-3mu[}}%
  \newcommand{\rrbracket}{\mathclose{]\mkern-3mu]}}}
\usepackage[hidelinks]{hyperref}
\newtheorem{theorem}{Theorem}
\newtheorem{proposition}[theorem]{Proposition}
\newtheorem{lemma}[theorem]{Lemma}
\newtheorem{corollary}[theorem]{Corollary}

\newcommand{\R}{\mathbb R}
\newcommand{\cH}{\mathcal H}

\newcommand{\cE}{\mathcal E}

\newcommand{\spt}{\operatorname{spt}}

\title{Sheetwise stationarity for stationary unions of $C^{1,\alpha}$ graphs}
\author{Mattia Luchese}
\date{}
\begin{document}
\maketitle

\begin{abstract}
We prove that stationarity of a finite union of $C^{1,\alpha}$ graphs implies stationarity of each individual graph, in arbitrary dimension and codimension.
The proof closely follows the approach of Minter and Stanbury \cite{MinterStanbury}, who also raised this question in their Remark~1.6.
As in their paper, we give a proof that uses the planar frequency results of Krummel, Minter and Wickramasekera \cite{KMW26} (but avoids the centre manifold construction of that paper). We also rely on the Schauder estimates for $Q$-valued functions established in \cite{Luc26}.
\end{abstract}

\paragraph{Introduction.}
This short note answers the question raised by Minter and Stanbury in \cite[Remark~1.6]{MinterStanbury}.
The first version was considerably longer, as it also contained the Schauder estimates for stationary $Q$-valued graphs that arose in the course of the proof.
Proving these estimates was in fact the main difficulty in the argument.
I subsequently extended them to a broader class of $Q$-valued quasilinear elliptic systems and published them in a separate paper \cite{Luc26}.
The present note should therefore be viewed as an application of that work.

The results of \cite{Luc26} in fact also give a route to sheetwise stationarity alternative to the one presented in this note. We briefly sketch this alternative route in the next paragraph.

Suppose that $u=\sum_{\ell=1}^Q\llbracket u_\ell\rrbracket$ is given as a sum of single-valued $C^{1,\alpha}$ sheets and that its graph is stationary.
Note that the branch set $\mathcal B_u$, in the sense of \cite[Section~2.1]{Luc26}, is empty, but the sheets $u_\ell$ may nevertheless have the same first jet at a point $y$ without coinciding in any neighbourhood of $y$.
Such points form the touching set $K_{\mathrm{gen}}=K\setminus K^\circ$ defined below.
For each $y\in K_{\mathrm{gen}}$, the point $Y=(y,u_1(y))$ is a singular point of the graph varifold $V$ with a multiplicity-$Q$ tangent plane, and hence belongs to the varifold branch set $\mathcal B_V$ in the sense of Krummel, Minter and Wickramasekera \cite{KMW26}.
Taking this tangent plane as horizontal and rescaling about $Y$ at a sufficiently small scale gives a varifold in the class $\mathcal V_L$ introduced in the proof of Theorem~5.7 of \cite{Luc26}.
As verified there, this class satisfies the hypotheses of Theorem~A of \cite{KMW26}.
Applying that theorem and projecting onto the base gives $\dim_{\mathcal H}K_{\mathrm{gen}}\leq n-2$.
Induction on $Q$ gives minimality of each sheet off $K_{\mathrm{gen}}$, and the dimension bound makes this set removable for the minimal surface system, yielding the conclusion that each $u_\ell$ is a smooth solution of the minimal surface system.

The proof sketched above, however, uses far more machinery than is needed for the problem at hand.
Indeed, it invokes Theorem~A of \cite{KMW26}, which gives a dimension bound for general $C^{1,\alpha}$ $Q$-valued maps with stationary graphs, even when these maps do not split locally into single-valued $C^{1,\alpha}$ sheets.
To prove Theorem~A in this generality, Krummel, Minter and Wickramasekera use a centre manifold construction going back to Almgren \cite[Chapter~4]{Alm00}.
Although they simplify Almgren's argument substantially, the construction remains a rather lengthy and involved part of their proof.

In Theorem~\ref{thm:main} below we give another proof that takes advantage of the fact that $u$ is already given as a sum of single-valued functions.
We still use the Schauder estimate of \cite[Corollary~5.1]{Luc26} and the planar frequency theory of \cite[Section~3]{KMW26}, but we can avoid the centre manifold construction.
This is very much in the spirit of the paper of Minter and Stanbury \cite{MinterStanbury}.
In their setting too, the dimension bound for the touching set already followed from Theorem~A of \cite{KMW26}, and they gave a simpler proof by PDE methods that avoids the centre manifold construction.
We follow the same philosophy here and, in fact, much of their proof quite closely.

The main differences are that we allow classical singularities of any density, work in arbitrary codimension, use the Schauder estimate of \cite[Corollary~5.1]{Luc26} in place of the $\varepsilon$-regularity theorem for stable hypersurfaces in \cite[Theorem~A]{MW24}, and avoid the Garofalo--Lin frequency argument by proving smoothness of the individual sheets before establishing the dimension bound for the full touching set.
We also give alternative arguments for a few technical steps in their proof.

\paragraph{Acknowledgements.}
I carried out this work during my PhD at the University of Cambridge.
I am very grateful to my supervisor, Neshan Wickramasekera, for many helpful conversations about this problem.
This work was supported by the UK Engineering and Physical Sciences Research Council (EPSRC) Studentship EP/T517847/1 and by the Cambridge Trust.

\paragraph{Statement on the use of AI.}
In writing the exposition, I used ChatGPT only for language editing.
The mathematical arguments were developed without AI assistance in March and April 2026.

\paragraph{Main result.}
Given $u\in C^1(B_1^n,\R^k)$, let $\Gamma_u=\{(x,u(x)):x\in B_1^n\}$.
We denote by $\mathbf v(u)$ the multiplicity-one varifold associated with $\Gamma_u$, given by
\[
\mathbf v(u)(\Phi)
=\int_{\Gamma_u}\Phi(X,T_X\Gamma_u)\,d\mathcal H^n(X)
\qquad\text{for }\Phi\in C_c\bigl((B_1^n\times\R^k)\times G(n+k,n)\bigr).
\]
The goal is to prove the following theorem.

\begin{theorem}\label{thm:main}
Let $n\geq2$, $k,Q\geq1$, $0<\alpha<1$, and
$u_1,\ldots,u_Q\in C^{1,\alpha}_{\mathrm{loc}}(B_1^n,\R^k)$.
If
\[
 V\coloneqq\sum_{i=1}^Q\mathbf v(u_i)
\]
is stationary in $B_1^n\times\R^k$, then each $\mathbf v(u_i)$ is stationary.
In particular, every $u_i$ is a smooth solution of the minimal surface
system in $B_1^n$.
\end{theorem}

For $P\in\mathbb R^{k\times n}$, set
\[
g(P):=I_n+P^TP,
\qquad
J(P):=\sqrt{\det g(P)},
\qquad
A(P):=D_PJ(P)=J(P)P\,g(P)^{-1}.
\]
We denote the entries of $g(P)^{-1}$ by $g^{ij}(P)$ and write $A_i^a(P):=D_{P_i^a}J(P)$ for the components of $A(P)$.
Recall that a $C^1$ function is a weak solution to the minimal surface system if
\[
\int_\Omega A(Du):D\varphi\,dx=0
\qquad\text{for every }\varphi\in C_c^1(\Omega,\mathbb R^k).
\]
We also let $K\subset B^n_1$ be the relatively closed set defined by
\[
K\coloneqq\{x\in B_1^n:u_1(x)=\cdots=u_Q(x),\ Du_1(x)=\cdots=Du_Q(x)\}.
\]

Before giving the proof of Theorem \ref{thm:main}, we recall the following two results. The first is a simple removability lemma.

\begin{lemma}[Removability]\label{lem:removability}
Let $E$ be relatively closed in an open set $\Omega\subset\R^n$, and let
$u\in C^1(\Omega,\R^k)$ be smooth and solve the minimal surface system
on $\Omega\setminus E$. Then $u$ solves the minimal surface system throughout $\Omega$ if either
\begin{enumerate}[label=\textup{(\roman*)},nosep]
\item $\cH^{n-1}(E)=0$; or
\item $E$ is locally contained in a finite union of embedded $C^1$
      hypersurfaces.
\end{enumerate}
\end{lemma}

\begin{proof}
(i) follows from a simple capacity argument. Fix $\varphi\in C_c^1(\Omega,\R^k)$. As $\mathcal H^{n-1}(E\cap \spt \varphi) =0$, it has vanishing $1$-capacity in $\Omega$. As it is also compact, there exists a sequence of cutoffs $\eta_j$ such that $0\leq\eta_j\leq1$, $\eta_j=0$ near $E\cap\spt\varphi$, $\eta_j\to1$ almost everywhere on $\Omega$, and
$\int_\Omega|D\eta_j|\to0$. Testing the minimal surface system off $E$ gives
\[
 \int_\Omega \eta_j \,A(Du):D\varphi
       =-\int_\Omega A(Du):(\varphi\otimes D\eta_j)\longrightarrow0.
\]
Since $A(Du)$ is continuous and hence locally bounded, the left-hand side converges to $\int_\Omega A(Du):D\varphi$.

For \textup{(ii)}, we argue locally by induction on the number of hypersurfaces.
Fix $x_0\in E$ and choose an open neighbourhood $U\subset\Omega$ such that $E\cap U\subset H_1\cup\cdots\cup H_N$, where the $H_j$ are embedded $C^1$ hypersurfaces.
After relabelling, we may assume that $x_0\in H_1$.
Shrinking $U$ if necessary, we may therefore assume that $H_1$ is relatively closed in $U$.

If $N=1$, write $H:=H_1$.
As in the proof of \cite[Lemma~4.3]{Luc26}, we test on the two
sides of $H$ with cutoffs vanishing in a shrinking neighbourhood of $H$.
The $C^1$ regularity of $H$ and the continuity of $A(Du)$ allow us to pass to the limit in the cutoff terms, obtaining the corresponding normal flux integrals over $H$.
These cancel because the outward unit normals are opposite, proving the weak equation across $H$.

For $N>1$, consider the open set $U\setminus H_1$.
The exceptional set $E\cap(U\setminus H_1)$ is relatively closed there and is contained in the remaining $N-1$ hypersurfaces.
The induction hypothesis therefore shows that $u$ solves the minimal surface system throughout $U\setminus H_1$, and elliptic regularity gives smoothness there.
The single-hypersurface case now extends the equation across $H_1$, giving the equation throughout $U$.
Since $x_0\in E$ was arbitrary, the equation holds throughout $\Omega$.
\end{proof}

The second result we need is a consequence of the planar frequency theory developed by Krummel, Minter and Wickramasekera in \cite{KMW26}.
It can be applied to our setting because of the Schauder estimate
established in \cite[Corollary~5.1]{Luc26}.
For $y\in K$, let $L_y(x):=u_1(y)+Du_1(y)(x-y)$ and, whenever $B_r^n(y)\Subset B_1^n$, set
\[
\cE(y,r)
:=
r^{-n-2}\sum_{\ell=1}^Q\int_{B_r^n(y)}|u_\ell-L_y|^2\,dx.
\]
Let also $K^\circ$ denote the interior of $K$ relative to $B_1^n$, consisting of the points that admit a neighbourhood on which $u_1=\cdots=u_Q$, and set $K_{\mathrm{gen}}:=K\setminus K^\circ$ (which is clearly relatively closed in $B_1^n$).
For every $y\in K_{\mathrm{gen}}$, the point $Y:=(y,u_1(y))$ is a density-$Q$ singular point of $V$ with a unique multiplicity-$Q$ planar tangent cone (indeed, regularity of $V$ near $Y$ would force all the sheets to coincide near $y$).

\begin{proposition}[Consequences of the planar frequency theory]
\label{prop:frequency-package}
Under the hypotheses of Theorem~\ref{thm:main}, there exists a finite-valued upper semicontinuous function
\[
\mathcal N:K_{\mathrm{gen}}\to[1,\infty)
\]
(called the ``planar frequency function'') with the following properties:
\begin{enumerate}[label=\textup{(\roman*)},nosep]
\item
\[
\dim_{\mathcal H}
\{y\in K_{\mathrm{gen}}:\mathcal N(y)\neq2\}
\leq n-2.
\]
\item
If the sheets are smooth and minimal near $y\in K_{\mathrm{gen}}$, then $\mathcal N(y)\geq2$.
\item
For every compact set $A\subset K_{\mathrm{gen}}$, there are $r_A>0$ and $C_A<\infty$, with $B_{r_A}^n(y)\Subset B_1^n$ for every $y\in A$, such that
\[
\cE(y,r)
\leq
C_A\left(\frac rR\right)^2\cE(y,R)
\]
whenever $y\in A$, $\mathcal N(y)\geq2$ and $0<r\leq R\leq r_A$.
\end{enumerate}
\end{proposition}

\begin{proof}
Fix $y\in K_{\mathrm{gen}}$.
After translating $(y,u_1(y))$ to the origin and rotating its common tangent plane to the horizontal plane, the varifold is locally represented by a $Q$-valued $C^{1,\alpha}$ map $w_y$, with $w_y(0)=Q\llbracket0\rrbracket$ and $Dw_y(0)=Q\llbracket0\rrbracket$.
After decreasing $\alpha$ if necessary, we may assume that it lies
in the admissible range of \cite[Corollary~5.1]{Luc26}
for every multiplicity $1\leq q\leq Q$.
There exists $r_y>0$ such that, for every $0<r\leq r_y$, the rescaled map
\[
v_{y,r}(x):=\frac{w_y(rx)}r
\]
is defined on $B_2^n$ and has sufficiently small slope, height excess and tilt excess, with the mass bounds required in \cite[Theorem~2.2]{KMW26}.
By \cite[Corollary~5.1]{Luc26}, we have
\[
\|v_{y,r}\|_{C^{1,\alpha}(B_{1/2}^n)}
\leq
C\|v_{y,r}\|_{L^2(B_1^n)},
\]
where $C=C(n,k,Q,\alpha)$ is independent of $y$ and $r$.
These are the quantitative graphical estimates used in the planar frequency arguments through the $\varepsilon$-regularity hypothesis.
Existence, finiteness, the bound $\mathcal N\geq1$ and upper semicontinuity of $\mathcal N$ follow by the arguments of \cite[Theorem~3.2 and Corollaries~3.4--3.5]{KMW26}.

For \textup{(i)}, we follow the stratification argument of \cite[Proposition~3.7]{KMW26}, using the upper semicontinuity along blow-up sequences established in \cite[Corollary~3.6]{KMW26} (the bound in \cite{KMW26} concerns the corresponding subset of the graph and transfers to the base because the projection map is Lipschitz).

For \textup{(ii)}, we observe that if the sheets are smooth near $y$, their first-order Taylor remainders are $O(|x-y|^2)$.
The first inequality in \cite[Corollary~3.4(ii)]{KMW26} then gives $\mathcal N(y)\geq2$.

For \textup{(iii)}, since $A$ is compact and each $u_\ell$ is $C^{1,\alpha}$, it is possible to choose the $r_y$ above to be equal to a common $r_A>0$ for all $y\in A$.
Moreover, finiteness and upper semicontinuity give a uniform upper bound for $\mathcal N$ on $A$.
Thus the second inequality in \cite[Corollary~3.4(ii)]{KMW26}, after returning to the original graphical coordinates, gives the conclusion whenever $\mathcal N(y)\geq2$ (using the compactness of $A$ again to find a uniform $C_A$ that works for all $y\in A$).
\end{proof}

\begin{proof}[Proof of Theorem~\ref{thm:main}]
We argue by induction on $Q$.
The case $Q=1$ follows from classical elliptic regularity.
Assume that $Q\geq2$ and that the theorem holds for all multiplicities smaller than $Q$.

\medskip
\noindent\emph{Step 1: Minimality outside $K_{\mathrm{gen}}$.}
Fix $x_0\notin K$ and partition the labels according to equality of their first jets at $x_0$.
Each group contains fewer than $Q$ labels.
Consider two labels $\ell$ and $m$ belonging to different groups.
If $u_\ell(x_0)\neq u_m(x_0)$, their graphs remain separated nearby.
Otherwise, $Du_\ell(x_0)\neq Du_m(x_0)$, so the implicit function theorem, applied to a suitable scalar component of $u_\ell-u_m$, shows that their coincidence set is locally contained in a $C^1$ hypersurface.
After shrinking to a common neighbourhood $U$ of $x_0$, all intersections between different groups are therefore contained in a finite union $E$ of relatively closed $C^1$ hypersurfaces in $U$.

On $U\setminus E$, the groups are locally separated in the ambient space, so stationarity of their union implies stationarity of each group.
The induction hypothesis then shows that each $u_\ell$ is smooth and solves the minimal surface system on $U\setminus E$.
Applying Lemma~\ref{lem:removability}\textup{(ii)} and elliptic regularity gives the same conclusion throughout $U$.
Thus every sheet is smooth and minimal on $B_1^n\setminus K$.

Near every point of $K^\circ$, all sheets coincide and $V=Q\mathbf v(u_1)$.
Stationarity therefore gives smoothness and minimality there as well.
Consequently, each $u_\ell$ is smooth and minimal on $B_1^n\setminus K_{\mathrm{gen}}$.

\medskip
\noindent\emph{Step 2: No subquadratic touching.}
Suppose that $\mathcal N(y_0)<2$ for some $y_0\in K_{\mathrm{gen}}$.
By upper semicontinuity, there is a ball $B_\rho^n(y_0)\Subset B_1^n$ such that $\mathcal N(y)<2$ for every $y\in K_{\mathrm{gen}}\cap B_\rho^n(y_0)$.
Proposition \ref{prop:frequency-package}\textup{(i)} then gives
\[
\mathcal H^{n-1}\bigl(K_{\mathrm{gen}}\cap B_\rho^n(y_0)\bigr)=0.
\]
By Step~1, each sheet is smooth and minimal outside this relatively closed set.
Lemma~\ref{lem:removability}\textup{(i)} and elliptic regularity therefore show that every sheet is smooth and minimal throughout $B_\rho^n(y_0)$.
Proposition~\ref{prop:frequency-package}\textup{(ii)} now gives $\mathcal N(y_0)\geq2$, a contradiction.
Hence $\mathcal N(y)\geq2$ for every $y\in K_{\mathrm{gen}}$.

\medskip
\noindent\emph{Step 3: Hessian bounds outside $K_{\mathrm{gen}}$.}
Fix $B_{2R}^n(x_0)\Subset B_1^n$ and set
\[
A:=K_{\mathrm{gen}}\cap\overline{B_{2R}^n(x_0)}.
\]
If $A=\varnothing$, the Hessians are bounded on $B_R^n(x_0)$ by the smoothness established in Step~1.
We therefore assume that $A\neq\varnothing$.
All constants below may depend on the chosen ball and the given sheets.

By Step~2 and Proposition \ref{prop:frequency-package}\textup{(iii)}, fixing the outer radius gives constants $C<\infty$ and $r_0\in(0,R)$ such that
\begin{equation}\label{eq:endpoint-L2-control}
\int_{B_r^n(z)}|u_\ell-L_z|^2\,dx
\leq Cr^{n+4}
\qquad
(z\in A,\ 0<r\leq r_0,\ \ell=1,\ldots,Q).
\end{equation}
Fix also a label $\ell$.

We bound the Hessian away from $K_{\mathrm{gen}}$.
Let $x\in B_{R}^n(x_0)\setminus K_{\mathrm{gen}}$ and set $d:=\operatorname{dist}(x,A)$.
If $d<r_0/2$, choose $z\in A$ with $|x-z|=d$.
The ball $B_{d/2}^n(x)$ lies in $B_{2R}^n(x_0)\setminus K_{\mathrm{gen}}$, where $u_\ell$ is smooth and solves the minimal surface system.
Hence, since $L_z$ is affine, each component $(u_\ell-L_z)^a$ of $u_\ell-L_z$ satisfies
\[
\sum_{i,j=1}^n g^{ij}(Du_\ell)\,D_{ij}(u_\ell-L_z)^a=0.
\]
The coefficients are uniformly elliptic and have uniformly bounded $C^{0,\alpha}$ norms by the initial $C^{1,\alpha}$ regularity.
Scalar interior estimates therefore give
\[
|D^2u_\ell(x)|
\leq
C d^{-2-n/2}
\|u_\ell-L_z\|_{L^2(B_{d/2}^n(x))}
\leq C,
\]
where the last inequality follows from \eqref{eq:endpoint-L2-control} and $B_{d/2}^n(x)\subset B_{2d}^n(z)$.
For points with $d\geq r_0/2$, the same bound follows from smoothness on a compact subset of $B_{2R}^n(x_0)\setminus K_{\mathrm{gen}}$.
Thus
\begin{equation}\label{eq:endpoint-hessian-bound}
|D^2u_\ell|\leq C
\qquad\text{on }B_{R}^n(x_0)\setminus K_{\mathrm{gen}}.
\end{equation}

\medskip
\noindent\emph{Step 4: Individual stationarity.}
Set
\[
\overline F:=\frac1Q\sum_{m=1}^Q A(Du_m),\qquad F_\ell:=A(Du_\ell)-\overline F.
\]
Each $F_\ell$ is continuous, vanishes on $K_{\mathrm{gen}}$ and is smooth on its complement.
By \eqref{eq:endpoint-hessian-bound} and the chain rule, $\int_{U\setminus K_{\mathrm{gen}}}|DF_\ell|^2\,dx<\infty$ for every $U\Subset B_1^n$.
It now follows from a classical zero extension lemma (recalled, for example, in \cite[Lemma~3.5]{Luc26}) that $F_\ell\in W^{1,2}_{\mathrm{loc}}(B_1^n,\mathbb R^{k\times n})$ and that $DF_\ell=0$ almost everywhere on $K_{\mathrm{gen}}$.
Since $\operatorname{div}F_\ell=0$ outside $K_{\mathrm{gen}}$ by Step~1, we conclude that $\operatorname{div}F_\ell=0$ distributionally throughout $B_1^n$.

For $\varphi\in C_c^1(B_1^n,\mathbb R^k)$, choose $\eta\in C_c^1(\mathbb R^k)$ equal to one near $\bigcup_{m=1}^Q u_m(\operatorname{spt}\varphi)$.
Testing stationarity of $V$ with $X(x,t):=(0,\eta(t)\varphi(x))$ gives $Q\int_{B_1^n}\overline F:D\varphi\,dx=0$, and hence $\operatorname{div}\overline F=0$.
Consequently, $\operatorname{div}[A(Du_\ell)]=\operatorname{div}F_\ell+\operatorname{div}\overline F=0$ in $\mathcal D'(B_1^n)$.
Every sheet therefore solves the minimal surface system weakly throughout $B_1^n$.
Elliptic regularity gives smoothness, and the corresponding graphs are individually stationary.
This completes the induction.
\end{proof}

\paragraph{Rectifiability of the touching set.}

We conclude the note by showing that the set $K_\mathrm{gen}$ is countably $(n-2)$-rectifiable. We need the following preliminary lemma, whose proof follows the argument of \cite[Lemma~1.9]{HS89}, which we recall here for completeness.

\begin{lemma}[Tangential contacts]
\label{lem:tangential-contact-capacity}
Let $\Omega\subset\R^n$ be connected and open, and let
$u,v\in C^\infty(\Omega;\R^k)$ be distinct solutions of the
minimal surface system. Then
\[
 C:=\{x\in\Omega:u(x)=v(x),\ Du(x)=Dv(x)\}
\]
is contained in the union of countably many smooth submanifolds of dimension $n-2$.
In particular, $C$ is countably $(n-2)$-rectifiable and has locally vanishing Sobolev $2$-capacity.
\end{lemma}

\begin{proof}
By the analyticity of smooth solutions of the minimal surface system \cite{Mor58}, the map $w:=u-v$ is real analytic.
Subtracting the equations in non-divergence form gives
\[
\sum_{i,j=1}^n a^{ij}D_{ij}w+\sum_{i=1}^n B^iD_iw=0,
\]
where $a^{ij}=g^{ij}(Du)$, the principal part is locally uniformly elliptic, and the $B^i$ are smooth $k\times k$ matrix-valued functions.
Since $w\not\equiv0$ and $\Omega$ is connected, every point of $C$ has finite vanishing order.
Write
\[
C=\bigcup_{m\geq2}C_m,
\qquad
C_m:=\{p:D^\gamma w(p)=0\text{ for }|\gamma|<m,\
D^mw(p)\neq0\}.
\]
Let $P_m$ be the nonzero leading homogeneous Taylor polynomial at $p\in C_m$.
Setting $w_\rho(x):=\rho^{-m}w(p+\rho x)$, the rescaled equation reads
\[
\sum_{i,j=1}^n a^{ij}(p+\rho x)D_{ij}w_\rho
+\rho\sum_{i=1}^n B^i(p+\rho x)D_iw_\rho=0.
\]
Since $w_\rho\to P_m$ in $C^2_{\mathrm{loc}}$ as $\rho\downarrow0$, passing to the limit gives
\begin{equation}\label{eq:leading-polynomial-equation}
\sum_{i,j=1}^n a^{ij}(p)D_{ij}P_m=0.
\end{equation}
Choose a component $a$ for which $P_m^a\not\equiv0$.
The vectors
\[
\{\nabla D^\beta w^a(p):|\beta|=m-1\}
\]
span a space of dimension at least two.
Indeed, otherwise, after rotating coordinates, they would all be parallel to $e_1$.
Consequently, $D_jD^\beta w^a(p)=0$ for every $j\geq2$ and $|\beta|=m-1$.
Since derivatives commute, every derivative of order $m$ involving a coordinate other than the first would vanish.
Thus $P_m^a(z)=c z_1^m$ for some $c\neq0$, contradicting \eqref{eq:leading-polynomial-equation} since $a^{11}(p)>0$ and $m\geq2$.

Choose two independent such vectors.
The implicit function theorem shows that the common zero set of the corresponding scalar derivatives is locally a smooth submanifold of codimension two containing $C_m$ in a neighbourhood of $p$.
A countable subcover for each $m$ proves that $C$ is contained in a countable union of smooth submanifolds of dimension $n-2$.

When $n\geq3$, the capacity conclusion follows from $\sigma$-finiteness of $\mathcal H^{n-2}$, \cite[Theorem~4.16]{EG15}, and countable subadditivity.
When $n=2$, it follows from countability and the zero $2$-capacity of a point (which can be checked using the logarithmic cutoff trick).
\end{proof}

\begin{corollary}
Under the hypotheses of Theorem~\ref{thm:main}, the set $K_{\mathrm{gen}}$ is countably $(n-2)$-rectifiable and has locally vanishing Sobolev $2$-capacity.
\end{corollary}

\begin{proof}
By Theorem~\ref{thm:main}, all sheets are smooth solutions of the minimal surface system.
If they all coincide identically, then $K_{\mathrm{gen}}=\varnothing$.
Otherwise, choose two distinct sheets $u_\ell$ and $u_m$.
Since $K_{\mathrm{gen}} \subset \{u_\ell=u_m,\ Du_\ell=Du_m\}$, the conclusion follows from Lemma~\ref{lem:tangential-contact-capacity}.
\end{proof}

\bigskip

\noindent
Department of Pure Mathematics and Mathematical Statistics, University of Cambridge\\
\textit{Email}: \texttt{ml2003@cam.ac.uk}


\begin{thebibliography}{KMW26}
\bibitem[Alm00]{Alm00}
F.~J.~Almgren, Jr.,
\emph{Almgren's Big Regularity Paper: $Q$-valued functions minimizing Dirichlet's integral and the regularity of area-minimizing rectifiable currents up to codimension $2$},
World Scientific Monograph Series in Mathematics, vol.~1,
World Scientific, River Edge, NJ, 2000.

\bibitem[EG15]{EG15}
L.~C.~Evans and R.~F.~Gariepy,
\emph{Measure Theory and Fine Properties of Functions},
revised edition, CRC Press, Boca Raton, FL, 2015.

\bibitem[HS89]{HS89}
R.~Hardt and L.~Simon,
\emph{Nodal sets for solutions of elliptic equations},
J. Differential Geom. \textbf{30} (1989), no.~2, 505--522.

\bibitem[KMW26]{KMW26}
B.~Krummel, P.~Minter and N.~Wickramasekera,
\emph{A branch set stratification for stationary varifolds with
epsilon-regularity}, preprint, 2026,
\href{https://arxiv.org/abs/2606.01511}{arXiv:2606.01511}.

\bibitem[Luc26]{Luc26}
M.~Luchese,
\emph{Perturbative Schauder estimates for $Q$-valued quasilinear elliptic
systems and a sharp dimension bound for branch sets of stationary graphs},
preprint, 2026,
\href{https://arxiv.org/abs/2609.12950}{arXiv:2609.12950}.

\bibitem[Mor58]{Mor58}
C.~B.~Morrey, Jr.,
\emph{On the analyticity of the solutions of analytic
non-linear elliptic systems of partial differential equations:
Part I. Analyticity in the interior},
Amer. J. Math. \textbf{80} (1958), no.~1, 198--218.

\bibitem[MS26]{MinterStanbury}
P.~Minter and S.~Stanbury,
\emph{Area minimising hypersurfaces mod $p$ do not admit immersed branch
points}, preprint, 2026,
\href{https://arxiv.org/abs/2603.03100}{arXiv:2603.03100}.

\bibitem[MW24]{MW24}
P.~Minter and N.~Wickramasekera,
\emph{A structure theory for stable codimension $1$ integral varifolds with applications to area minimising hypersurfaces mod $p$},
J. Amer. Math. Soc. \textbf{37} (2024), no.~3, 861--927.

\end{thebibliography}
\end{document}